\documentclass[12pt]{amsart}

\usepackage{cite}

\usepackage{kotex}
\usepackage{lineno,hyperref}
\modulolinenumbers[5]
\usepackage{epstopdf}
\usepackage{pifont}
\usepackage{latexsym}
\usepackage{graphics}
\usepackage{graphicx}
\usepackage{float}
\usepackage[dvips]{xy}
\xyoption{all}
\usepackage{ifpdf}
\usepackage{multirow}
\usepackage{amssymb, amsthm, amsmath}
\usepackage{hhline}
\usepackage{centernot}
\usepackage{verbatim, comment, color}
\usepackage{enumerate}
\usepackage{setspace}
\usepackage{lipsum}
\usepackage{subcaption}

\usepackage{mathtools}

\newtheorem{theorem}{Theorem} 
\newtheorem{proposition}[theorem]{Proposition}
\newtheorem{lemma}[theorem]{Lemma}
\newtheorem{remark}[theorem]{Remark}
\newtheorem{corollary}[theorem]{Corollary}
\newtheorem{definition}[theorem]{Definition}

\newcommand{\ba}{\begin{align}}
\newcommand{\ea}{\end{align}}  

\newcommand{\be}{\begin{equation}}
\newcommand{\dd}{\rho}
\newcommand{\ee}{\end{equation}}
\newcommand{\bea}{\begin{eqnarray}}
\newcommand{\eea}{\end{eqnarray}}
\newcommand{\barr}{\begin{array}}
\newcommand{\earr}{\end{array}}
\newcommand{\bn}{\begin{enumerate}}
\newcommand{\en}{\end{enumerate}}
\newcommand{\bi}{\begin{itemize}}
\newcommand{\ei}{\end{itemize}}
\newcommand{\bbbm}{\begin{pmatrix}}
\newcommand{\eeem}{\end{pmatrix}}

\newcommand{\bbN}{{\bf N}}

\newcommand{\bbP}{{\bf P}}

\newcommand{\bbR}{{\bf R}}
\newcommand{\bbS}{{\bf S}}

\newcommand{\cP}{{\cal P}}

\newcommand{\cPo}{\cP_{{\rm on}}}
\newcommand{\cPb}{\cP_{{\Delta}}}

\newcommand{\cM}{{\cal M}}

\newcommand{\R}{{\mathbf R}}

\newcommand{\al}{\alpha}

\newcommand{\bt}{\beta}

\newcommand{\De}{\Delta}
\newcommand{\ald}{{\al_{\De^d}}}
\newcommand{\ep}{\epsilon}

\newcommand{\la}{\lambda}
\newcommand{\La}{\Lambda}

\newcommand{\ignore}[1]{}{}
\newcommand{\noin}{\noindent}

\newcommand{\nn}{\nonumber}

\newcommand{\h}{\bar}

\newcommand{\q}{\quad}

\newcommand{{\QED}}{{\hfill QED} \smallskip}

\newcommand{\Rn}{{\R^n}}
\newcommand{\spt}{\mathop{\rm spt}}
\newcommand{\Tr}{\mathop{\rm Tr}}

\renewcommand{\subset}{\subseteq}

\renewcommand{\phi}{\varphi}
\newcommand{\cal}{\mathcal}

  \DeclareMathOperator*{\argmax}{argmax}

 \definecolor{darkspringgreen}{rgb}{0.09, 0.45, 0.27} 
 \definecolor{darkgray}{rgb}{0.66, 0.66, 0.66}
\numberwithin{equation}{section}
\numberwithin{theorem}{section}
    
\begin{document}
\title
[Maximizing powers of the angle between projective points]
{Maximizing expected powers of the angle between pairs of points in projective space$^*$} 

\thanks{$^*$With an appendix by Dmitriy Bilyk, Alexey Glazyrin, Ryan Matzke, Josiah Park and Oleksandr Vlasiuk.\\
\em TL is grateful for the support of ShanghaiTech University, and in addition, to the University of Toronto and its Fields Institute for the Mathematical
Sciences, where parts of this work were performed.  RM  acknowledges partial support of his research by the Canada Research Chairs Program and
Natural Sciences and Engineering Research Council of Canada Grants RGPIN 2015-04383 and 2020-04162. 
\copyright 2021 by the authors.
}
\date{\today}

\author{Tongseok Lim and Robert J. McCann}
\address{Tongseok Lim: Krannert School of Management \newline  Purdue University, West Lafayette, Indiana 47907}
\email{lim336@purdue.edu}
\address{Robert J. McCann: Department of Mathematics \newline University of Toronto, Toronto ON Canada}
\email{mccann@math.toronto.edu}

\begin{abstract}
Among probability measures on $d$-dimensional real projective space,  one which maximizes
the expected angle $\arccos(\frac{x}{|x|}\cdot \frac{y}{|y|})$ between independently drawn projective points $x$ and $y$
was conjectured to equidistribute its mass over the standard Euclidean basis $\{e_0,e_1,\ldots, e_d\}$ 
 by Fejes T\'oth \cite{FT59}.
 If true,  this conjecture evidently implies the same measure maximizes the expectation of 
$\arccos^\alpha(\frac{x}{|x|}\cdot \frac{y}{|y|})$ for any exponent $\alpha > 1$.    
The kernel $\arccos^\alpha(\frac{x}{|x|}\cdot \frac{y}{|y|})$ represents the objective of an infinite-dimensional quadratic program.
We verify discrete and continuous versions of this {milder} conjecture in a non-empty range $\alpha > \alpha_{\Delta^d} \ge 1$,  
and establish uniqueness of the resulting maximizer $\hat \mu$ up to rotation. 
{We show $\hat \mu$ no longer maximizes when $\alpha<\alpha_{\Delta^d}$.}
At the endpoint $\alpha=\alpha_{\Delta^d}$ of this range, we show 
another maximizer $\mu$ must also exist which is not a rotation of $\hat \mu$.
For the continuous version of the conjecture, 
{an appendix} provided by Bilyk et al in response to an earlier draft of this work 
{combines with the present improvements
to yield} $\alpha_{\Delta^d}<2$.
The original conjecture $\ald=1$ remains open (unless $d=1$).

{However, in the maximum possible range $\alpha>1$,} we show $\hat \mu$ and its rotations 
maximize the aforementioned expectation uniquely on a sufficiently small ball
in the $L^\infty$-Kantorovich-Rubinstein-Wasserstein metric
$d_\infty$ from optimal transportation;  the same is true for
any measure $\mu$ which is mutually absolutely continuous with respect to $\hat \mu$, but the size of the
ball depends on {$\alpha,d$, and} $\|\frac{d \hat \mu}{d\mu}\|_{\infty}$. 
 \end{abstract}

\maketitle


\section{Introduction}

Our previous work explored the embedding of various conjectures by L.~Fejes T\'oth into one-parameter families of maximization problems parameterized by an exponent $\al \in \R$ \cite{LM20_1}.  As the parameter is varied,  the maximizer in each family might theoretically bifurcate continuously or discontinuously, as caricatured by the following familiar examples:
for $f_\la(x) = -\frac14x^4 - \frac\la 2 x^2$ the maximizer $x_\la \in \argmax_\R f_\la$ bifurcates continuously at
$\la=0$ (in a so-called pitchfork), whereas for $f_\la(x) = -\frac14(1-  x^2)^2 - \la x$ the maximizer changes discontinuously from $x_\la \sim -1$ to $x_\la \sim +1$ as $\la$ increases through $0$.  In the current manuscript,  we establish existence and finiteness
of a critical threshold
for the exponent $\al$ in each family, at which 
our maximizer changes discontinuously unless the corresponding Fejes T\'oth conjecture is true.  An appendix by Bilyk, Glazyrin, Matzke, Park, and Vlasiuk leads to an effective bound  for the critical
threshold in the most symmetrical of these families. 

 
Choose $N$ unoriented lines through the origin of $\R^{d+1}$.  
 Suppose each pair of lines repel each other with a force {whose strength is} 
 independent of the (acute) angle
between them,  so that they prefer to be orthogonal to each other.  However,  unless $N \le d+1$,
it is impossible for all pairs of lines to be orthogonal.  What then are their stable configurations?
This question is partly 
addressed by an unsolved conjecture, which for $d=2$ dates back to Fejes T\'oth \cite{FT59}:
namely that the sum of the angles between these lines is maximized if the lines are distributed as evenly as possible amongst the coordinate axes of some orthonormal
basis for $\R^{d+1}$.  
For $d\ge 2$ this conjecture has motivated 
a recent series of works by  Bilyk, Dai, Glazyrin, Matzke, Park, Vlasiuk in different combinations
 \cite{BD19} \cite{BGMV19} \cite{BGMV19-2} \cite{BM19}, and by Fodor, V\'{\i}gh and Zarn\'{o}cz \cite{FVZ16}.
Other authors have also considered versions of the problem
for oriented as well as unoriented lines, in the limit $N=\infty$  --- where it can be formulated as an infinite-dimensional quadratic program ---
 and/or with different powers $\al$
of the angle or distance between them, e.g. \cite{PS31} \cite{Bjorck56} \cite{AS74} \cite{BDM18}.
In recent work we verified the unoriented conjecture for bounded and unbounded values of $N$ in the mild repulsion limit $\al=\infty$ \cite{LM20_1}. In the present manuscript we shall show there are threshold exponents $\ald(N) \in [0,\infty)$ and $\ald=\ald(\infty) \in [1,\infty)$, such that 
the same configurations continue to be optimal  precisely when $\al \ge \ald{(N)}$,
and uniquely optimal apart from known symmetries if and only if $\al>\ald{(N)}$.  
An appendix provided by Bilyk et al 
in response to an earlier draft of this work {combines with the present improvements to imply} $ \ald(N) = \ald<2$
for all $d \in \bbN$ and $N$ divisible by $d+1$, 
but the conjectured values $\ald=1=\ald(N)$ for all $N > d+1$ remain open unless $d=1$. However,
we show the same configuration with arbitrary positive weights 
remains {\em locally} optimal --- and hence stable --- in a suitable sense over the maximal range $\al > 1$.

Such explicit characterizations of optimizers which display symmetry breaking are rare in the context of infinite-dimensional
quadratic progamming,  and potentially valuable.  In modern data sciences it has also become more and more important
to study optimization problems in the curved geometries which arise from dimensional reduction,  
and we hope our approach to this model problem
may also prove relevant in such settings.

The Euclidean unit sphere $\bbS^d= \{x \in \R^{d+1} \ | \ |x|=1 \}$ 
gives a double cover for the real projective space $\R\bbP^d:= \bbS^d / \{\pm\}$
with covering map $x \to \tilde x := \{x,-x\}$.  Let $\rho$ be the geodesic distance 
$\dd(x,y) = \arccos (x \cdot y)$ on $\bbS^d$.
Define $\La_0 : [0, \pi] \to \R$ and $\La : \bbS^d \times \bbS^d \to \R$ by
\begin{align}\label{Lambda0}
   & \La_0(t):= \frac2\pi \min\{t,\pi -t\} \\ 
   &\La(x,y) := \La_0( \rho(x,y) ),
\label{projective distance}
\end{align}
so that $\La(x,y)$ is the distance on projective space,
rescaled to have unit diameter.  Note $\La(x,y)=1$ if and only
if $x\cdot y=0$.  

Let  $\cM(\bbS^d)$ denote the set of signed Borel measures on $\bbS^d$ with finite total variation and $\cP(\bbS^d) := \{ 0 \le \mu \in \cM(\bbS^d) \ | \ \mu(\bbS^d)=1\}$ the subset of probability measures.  
 Let $\spt(\mu)$ denote the smallest closed subset of $\bbS^d$ containing the full mass of $\mu$.
For $\al>0$ define the bilinear form
\be\label{quadratic}
B_{\al}(\mu, \nu) := \iint \La^\al(x,y) d\mu(x)d\nu(y), \q \mu, \nu \in \cM(\bbS^d),
\ee
and the corresponding energy
\be\label{energy}
E_{\al}(\mu) :=  \frac{1}{2} B_\al(\mu,\mu)
\q \al \in (0, \infty].
\ee
In particular, notice $\La^\infty(x,y)=1$ if $x \cdot y =0$, and zero otherwise.
\\

 For some $1 \le \ald <2$,  we shall establish
\begin{equation}\label{maximum value}
\max_{\mu \in \cP(\bbS^d)} E_\alpha(\mu) = \frac{d}{2d+2}
\ \mbox{\rm if and only if}\ \al \ge \ald, 
\end{equation}
and describe the set of maximizers precisely in the interior $\al > \ald$ 
of this range. 

Write  $\mu \equiv \nu$, and say the measures $\mu,\nu \in \cP(\bbS^d)$ are {\em essentially equivalent},
if and only if there is a rotation $M\in SO(d)$ such that $\mu(A \cup -A) = \nu(M(A \cup -A))$ for each open set
$A \subset \bbS^d$.  In other words,  $\mu \equiv \nu$ if and only if
their canonical projections 
onto $\R\bbP^d$ 
are rotations of each other.  When a maximum over measures is attained uniquely up to essential equivalence,
we say the maximizer is {\em essentially unique}.

 For $N \in \bbN :=\{1,2,\ldots\}$ consider the following collections of discrete probability measures
 on the $d$-sphere:
\begin{align*}
\cP_N^=(\bbS^d) &:= \{ \mu \in \cP(\bbS^d) \ | \  \mu=\frac{1}{N}\sum_{i=1}^N \delta_{x_i}, \ \ x_i \in \bbS^d \}, 
\nonumber \\
\cPo(\bbS^d) &:= \{ \mu \in \cP(\bbS^d) \ | \ \spt(\mu) \text{ is an orthonormal basis of } \R^{d+1} \}, 
\nonumber \\ 
\cPo^=(\bbS^d) &:= \{ \mu  \in \cPo(\bbS^d) \ | \ \mu[\{x\}] = \frac{1}{d+1} \ \forall\ x \in \spt \mu \},
\nonumber \\
\cPb(\bbS^d) &:= \{ \mu \in \cP(\bbS^d) \ |   \ \mu \equiv \nu \ \text{ for some } \ \nu \in \cPo(\bbS^d) \},
\\
\cPb^=(\bbS^d) &:= \{ \mu \in \cP(\bbS^d) \ | \ \mu \equiv \nu \ \text{ for some } \ \nu \in \cPo^=(\bbS^d) \}.
\end{align*}
The closure $\overline{\cPb(\bbS^d)}$ of $\cPb(\bbS^d)$ in the narrow topology denotes the set of measures whose projective support is contained in an orthonormal basis 
$V=\{v_0,\ldots,v_d\}$ --- which we think of as forming the vertices of a maximal projective simplex  --- and 
$\cPb^=(\bbS^d)$ is the subset of measures which equidistribute their mass over these vertices.
 
Motivated by \cite{LM20_1}, the following elementary lemma is our starting point:
 it establishes existence of a threshold for the transition that we intend to describe: 
 
\begin{lemma}[Threshold exponent for equidistribution over maximal simplices]\label{iff} 
There exists a unique $\ald \in [1,\infty]$ 
such that  $0<\al < \ald$ implies $\cPb^=(\bbS^d)$ disjoint from $\argmax_{\cP(\bbS^d)} E_\al$, whereas $\ald<\al$ implies
\be\label{unique max}
\cPb^=(\bbS^d) = \argmax_{\cP(\bbS^d)} E_\al.
\ee

Moreover,  \eqref{unique max} fails when $\al=1$, meaning $\cPb^=(\bbS^d)$ is either strictly contained in or disjoint from 
$\argmax_{\cP(\bbS^d)} E_1$.
\end{lemma}

\noin{\bf Proof.} Fix $0<\al_0<\infty$.  Notice $\Lambda(x,y) \le 1$ on $\bbS^d \times \bbS^d$,  with equality if and only if $x \cdot y=0$.  
For each $\mu \in \cP(\bbS^d)$, it follows that $\al \in (0,\infty] \mapsto E_\al(\mu)$ is nonincreasing; moreover, it
decreases strictly unless each distinct pair of points $x\ne y$ in $\spt \mu$ are orthogonal,
i.e.\ unless $\mu \in \overline{\cPb(\bbS^d)}$.  
On $\overline{\cPb(\bbS^d)}$ it is well-known \cite{LM20_1} and easy to check 
\begin{equation}\label{Perron-Frobenius}
\cPb^=(\bbS^d) = \argmax_{\mu \in \overline{\cPb(\bbS^d)}} E_{\al_0}(\mu)
\end{equation}
using, e.g. the Perron-Frobenius theorem.

Now suppose some (and hence all) $\hat \mu \in \cPb^=(\bbS^d)$ attain the maximum \eqref{maximum value} when $\al=\al_0$.
For $\al>\al_0$ and $\mu \in \cP(\bbS^d)$,  the monotonicity mentioned above implies 
$E_\al(\mu)\le E_{\al_0}(\mu) \le E_{\al_0}(\hat \mu) = E_\al(\hat \mu)$.  The first inequality is strict unless $\mu \in \overline{ \cPb(\bbS^d)}$;
for $\mu \in \overline{\cPb(\bbS^d)}$, the second inequality becomes strict unless $\mu \in \cPb^=(\bbS^d)$, according to \eqref{Perron-Frobenius}.
Thus we conclude \eqref{unique max} for all $\al>\al_0$.  Taking the infimum of such $\al_0$ yields $\ald$, with the usual convention $\ald=\infty$
if no such $\al_0$ exists.

It remains to show $\ald \ge 1$.  When $d=1$ the energy $E_1(\mu)$ takes the constant value $1/4$ on the closed convex hull of  $\cPb^=(\bbS^1)$ --- which consists of all measures on $\bbS^1/\{\pm\}$
that are invariant under a rotation by angle $\pi/2$ --- thus includes not only $\hat \mu$ but also the uniform measure on $\bbS^1$.  

For $d>1$,  it therefore follows that $E_1(\hat \mu) = E_1(\mu)$,  where
$\hat \mu = \frac1{d+1} \sum_{i=0}^d \delta_{v_i} \in \cPb^=(\bbS^d)$ and $\mu$ represents the average of $\hat \mu$ with any (or all) of its rotations in the plane spanned by 
$v_0$ and $v_1$, say. This shows \eqref{unique max} cannot hold when $\al=1$,  hence $\ald \ge 1$ as desired. \QED
\\

We call the assertion $\ald = 1$ (regardless of dimension $d$) the {\em continuous Fejes T\'oth conjecture}, which has remained open at least since 1959 \cite{FT59}, except for one dimension $d=1$.

A major goal of the present manuscript is to show 
that at this threshold exponent provided by Lemma \ref{iff} 
the inclusion 
\be\label{nonunique max}
\cPb^=(\bbS^d) \subset \argmax_{\cP(\bbS^d)} E_\ald 
\ee
becomes strict (Theorem \ref{main1}), meaning the optimizer is no longer essentially unique. For $\ald<\infty$, the inclusion follows from $\Gamma$-convergence ideas explained in Section \ref{S:nonunique}; however,
its strictness requires us to first establish a new local {stability} result: Theorem \ref{localmax}.
This local stability holds {in the maximum possible range $\al>1$, thus is} 
consistent with the conjectured value $\ald=1$; when combined  with the 
finiteness of $\ald$ 
shown subsequently 
 it implies discontinuous dependence of the {maximizer} on $\al$ at $\ald>1$ otherwise (Corollary \ref{C:bifurcate}).

When $\al=\infty$, the essential unique characterization \eqref{unique max} of the optimizer follows from our earlier works 
\cite{LM20_1} \cite{LM20_3} together with an estimate $N(d)$ on the cardinality of support of local energy maximizers
for $\al \ge 4$ which develops from the ideas of various authors \cite{CarrilloFigalliPatacchini17} \cite{KKS19} \cite{V20}.
The local stability result established hereafter then implies $\ald<\infty$.
Following our announcement of this result,
Bilyk et al reported to us a simple majorization argument that immediately gives $\ald \le 2$;
they have kindly agreed to include their argument in the present paper in the form of a separately 
authored appendix.
Motivated in part by this development,  we have refined our local stability analysis below to imply $\ald<2$.  
Indeed, for every $\al > 1$, our stability analysis
now shows each measure in the broader class $\cPb(\bbS^d)$ {of unbalanced simplices}
is an essentially unique local maximizer for the energy $E_\al(\mu)$
on an appropriately metrized version of the landscape $\cP(\bbS^d)$.  
To formulate both this stability result (Theorem \ref{localmax})
and the $\Gamma$-convergence statement for the family of energies $(E_\al)_{\al>0}$
requires us to recall the metrics $ \{d_p\}_{p \in [1,\infty]}$ on this landscape from the theory of optimal transportation which we use to define the relevant topologies. 
Unless the conjecture $\ald=1$ holds, we show (i) the $d_\infty$-continuity of any curve $(\mu_\al)_{\al >0}$ of optimizers
\[
\mu_\al \in \argmax_{\cP(\bbS^d)} E_\al
\]
must fail at $\al=\ald$ {(Corollary \ref{C:bifurcate})}, and (ii) when $\al=\ald$ the optimizers 
form a disconnected set in the $d_\infty$-topology {(Corollary \ref{C:disconnect})}.

These theorems {partly} echo our results concerning particles interacting via 
strongly attractive / mildly repulsive potentials on Euclidean space \cite{LimMcCann19p}. 
In the present context,  the interaction kernel $\Lambda^\al$ acts purely repulsively on $\R\bbP^d$ (or attractive-repulsively on $\bbS^d$), with compactness of the space substituting
for strong attraction at large distances.  The restriction $\al \ge 2$ corresponds to the mildly repulsive range of potentials from \cite{BalagueCarrilloLaurentRaoul13}  \cite{CarrilloFigalliPatacchini17} 
 \cite{LimMcCann19p}. {A remarkable difference between our results and those of the Euclidean setting, however, is that the characterization of the optimizers \eqref{unique max} continues to hold in the interval $\al \in (\ald,2) \ne \emptyset$.  This 
 helped motivate our effort to exclude the case of equality from the estimate  $\ald \le 2$ provided by  Bilyk et al's appendix.  By contrast,   for the analogous problem on Euclidean space,  in which the charges interact to minimize the attractive-repulsive pair potential $W(x) = \phi(x) - \frac1\al |x|^\al$,  Balagu\'e et al \cite{BalagueCarrilloLaurentRaoul13} have shown 
 that, e.g., $\phi \in C^{1,1}_{loc}(\Rn)$ and $\al >2-d$ imply any ($d_\infty$-local) optimizer with compact support has Hausdorff dimension at least $2 -\al$. 
 We sketch their argument briefly:  the Euler-Lagrange equation satisfied by an optimizer $\mu \in \cP(\R^d)$ yields $\Delta (W* \mu) \ge 0$ on $\spt \mu$.
Now
\[
\Delta \frac{|x|^{\al}}{\al{ (\al-d+2)}} = |x|^{\al-2} + {c_d(\al)} \delta_0 
\ge |x|^{\al-2}
\]
so
\begin{eqnarray*}
{(\al-d+2)} {\iint_{\R^d \times \R^d} |x-y|^{\al-2} d\mu(x) d\mu(y) }
\le  \|\Delta \phi\|_{L^\infty(\R^d)}.
 \end{eqnarray*}
 But {finiteness} of this {singular integral} is known to imply any set of full $\mu$ measure has Hausdorff dimension at least $2-\al$ (i.e.~positive for $\al < 2$).  
 Our results 
 show the analogous conclusion fails in a nonempty range $\al\in[\ald,2)$ on the projective sphere,  reflecting its compact geometry and providing some  further evidence supporting} the continuous {\em Fejes T\'oth conjecture} $\ald=1$.  

 Lastly, we note that Fejes T\'oth actually considered the maximization problem over the discrete domain $\cP_N^=(\bbS^d)$ for each 
$N \in \bbN$. In this setup it seems that majorization approach which Bilyk et al describe in the appendix below does not apply (unless $N$ is a multiple of $d+1$). In Section \ref{S:discrete} 
we shall address this by giving corresponding results regarding the existence of an $N$-dependent threshold exponent
$\ald (N) \in [1,\infty)$ for this problem when $N>d+1$  with Fejes T\'oth's conjectures becoming equivalent to the assertion
$\ald(N)=1$; see Lemma \ref{L:N-iff} and Theorem \ref{main2}.

\section{Measures supported on an orthonormal basis are local energy maximizers}
\label{S:local}

 In this section we show that the conjectured maximizer of \eqref{maximum value} for $\al=1$ is --- up to projective rotations
--- a strict local maximizer for all $\al>1$.  Moreover, the same is true for all measures which share its support.
For $1 \le p < +\infty$ 
 define the $L^p$-Kantorovich-Rubinstein-Wasserstein (optimal transport) distance between 
 $\mu,\nu \in \cP(\bbS^d)$ by
\be\label{KRW metric}
d_p(\mu,\nu) 
:=\inf_{\gamma \in \Gamma(\mu,\nu)} 
\left(\int_{\bbS^d \times \bbS^d} \rho(x,y)^p d\gamma(x,y)\right)^{1/p},
\end{equation}
where the infimum is taken over the set $\Gamma(\mu,\nu)$ of joint probability measures on $\bbS^d \times \bbS^d$
having $\mu$ and $\nu$ as their left and right marginals.
 For $p\ne \infty$, the distance $d_p$ is well-known to metrize
narrow convergence (against continuous bounded test functions,
e.g.~Theorem 7.12 of \cite{Villani03}),
so $\cP(\bbS^d)$ becomes a compact metric space under $d_p$.
The limit 
\be\label{sand metric}
d_\infty(\mu,\nu) 
:= \inf_{\gamma \in \Gamma(\mu,\nu)} \sup_{(x,y)\in \spt \gamma} \dd(x,y)
\end{equation}
is also a distance,  but metrizes a much finer and non-compact topology on $\cP(\bbS^d)$
  (but not on $\cP_N^=(\bbS^d)$, as we shall see in Lemma \ref{L:discrete narrow}).
For  $\al >1 $, this finer topology allows us  to establish local energy maximality
of even the wildly unbalanced measures in $\cPb(\bbS^d)$,
echoing its uses in other settings \cite{McCann06} \cite{BalagueCarrilloLaurentRaoul13} \cite{LimMcCann19p}.

   Let $e_0 = (0,0,...,0,1), e_1= (1,0,...,0), \dots, e_d=(0,...,0,1,0)$ be the standard basis of $\R^{d+1}$.    Let $D(x,r)$ and $B(x,r)$    denote the open balls of center $x$ and radius $r$ in $\bbS^d$ and $\R^d$ respectively,
and denote the Euclidean norm by $| \cdot |$. The following lemma is standard:

\begin{lemma}[Barycenter characterization]
\label{jensen}
Let $\mu$ be a Borel probability measure on $\R$ and $ c \in L^1(\R,d\mu)$. Then 
\be
 \int |x-c(y)|d\mu(y) \ge\Big|x-\int c(y) d\mu(y)\Big| \ \ \text{for all} \ \ x \in \R.
\ee
\end{lemma}

 \noin{\bf Proof}: Jensen's inequality.
\QED
\\

The function $x\in \bbS^d \mapsto \Lambda(x,e_1)$ peaks like a roof at the equator $e_1^\perp \cap \bbS^d$, where $x^\perp:= \{ y \in \R^{d+1} \mid x \cdot y=0\}$.
Thus, near $e_0 \in \bbS^d$, the potential corresponding to unit Dirac masses at each of the other basis vectors $e_i$
behaves asymptotically like the $\ell^1$ norm on the tangent space $T_{e_0}\bbS^d$,
 which obviously dominates the $\ell^2$ (Riemannian) norm there.
Thus $C'<1$ yields
\be\label{ell1 model}
\frac\pi 2 \sum_{i=1}^d  (1-\La(x,e_i)) \sim \sum_{i=1}^d | x \cdot e_i| 
\ge C' \Big(\sum_{i=1}^d | x \cdot e_i|^2\Big)^{1/2},
\ee
with the final inequality being strict unless both sides vanish,  in which case $x$ is a multiple of $e_0$.
The next lemma gives a perturbation of this result that plays a crucial role in our local stability argument.

\begin{proposition}[Recentered potentials near conjectured optimizer]
\label{P:aggregation}
Given $0<C<\frac{2}{\pi}$, there exists $r(d,C)>0$  such that: 
if $\nu_i \in \cP(D(e_i,r))$ for $r<r(d,C)$ and each $i=1,...,d$,
(i.e. $\nu_i$ is a probability measure on $\bbS^d$ with $\nu_i(D(e_i,r))=1$),
then there exists $\bar x = \bar x(\nu_1,\ldots,\nu_d) \in \bbS^d$ satisfying
\be\label{agg}
\sum_{i=1}^d \int (1 - \La(x,y) ) d\nu_i(y) \ge C\dd(x,\bar x) \q \forall \, x \in D(e_0,r).
\ee
Finally, $\bar x \to e_0$ as $r \to 0$.
\end{proposition}

\noin{\bf Proof.} 
Fix the Riemannian exponential map 
\begin{equation}\label{spherical exponential}
\exp_{e_0} v = e_0 \cos |v| + \frac{v}{|v|} \sin |v| 
\end{equation}
from the tangent space $T_{e_0}M = e_0^\perp \cong \bbR^d \subset \bbR^{d+1}$  onto $\bbS^d$.
Let $i \in \{1,...,d\}$ and $r$ be small. For each $y \in D(e_i,r)$ and small $v \in e_0^\perp$, consider the function $V_y(v) := 1 - \La(\exp_{e_0} v,y)$ and its zero set 
$Z_y \subset e_0^\perp$.  Note $ {\exp_{e_0} Z_y =  y^\perp\cap \bbS^d} = \{x  \ | \ \La(x,y) = 1 \} = \{ x  \ | \ \dd(x,y) = \frac{\pi}{2} \}$.

Let $Q$ be the hypercube $[-s,s]^d \subset e_0^\perp$ for $s >0$. For $r \ll s \ll 1$ both sufficiently small,  $Z_y \cap Q$ can be viewed as the graph of a smooth function $\zeta_y : [-s,s]^{d-1} \to [-s,s]$, that is, $v=(v_1,...,v_d) \in Z_y$ if and only if $ v_i=\zeta_y(v_{-i})$ where $ v_{-i} := (v_1,...,v_{i-1},v_{i+1},...,v_d)$. 
 Here the $\zeta_y$ are uniformly smooth, hence $\zeta_y  \to 0 = \zeta_{e_i}$ in smooth norms as $r  \to 0$,  
 our notation $r \ll s$ signifying that $r$ tends to zero faster than $s$.
Since the linearization of \eqref{spherical exponential} at $v=0$ gives the identity map, from the definition of $\La$ and $C < \frac{2}{\pi}$, observe if $r,s$ are sufficiently small
\be
V_y( v) \ge C|f_y(v)| \q \text{for all} \   v \in  Q \nn
\ee
where $f_y (v) := v_i - \zeta_y(v_{-i})$ on $v=(v_1,\ldots,v_d) \in Q$.
This implies
\be
F_i(v) := \int V_y (v) \, d\nu_i(y) \ge C \int |f_y (v) |\, d\nu_i (y) \q  \text{for all} \   v \in Q. \nn
\ee
Let $\xi_i :  [-s,s]^{d-1} \to [-s,s]$ be defined by $\xi_i(v_{-i})= \int \zeta_y(v_{-i}) d\nu_i(y) $, and $Z^i \subset Q$ be its graph $Z^i:= \{(v_1,...,v_d) \in Q \mid v_i=\xi_i( v_{-i} )\}$. Define  a function $g_i$ on $Q$ by $ g_i(v_1,\ldots,v_d) = v_i - \xi_i(v_{-i})$. 
Then by Lemma \ref{jensen}, we find
\be
\int |f_y| \, d\nu_i (y) \ge |g_i| \q \text{on} \ Q. \nn
\ee
 From uniformity in $i$ of smoothness of $\xi_i$,  it follows that $|\xi_i|$ and $|\nabla \xi_i|$ 
converge uniformly to zero as $r \to 0$ (regardless of the choice of $\nu_i$), and their graphs are orthogonal in the limit.
Thus $v^r := \bigcap_{i=1}^d Z^i $ is a singleton for  small enough $r$ and $v^r \to 0$ as $r \to 0$.

Composing the map $v \in e_0^\perp \mapsto G^r(v)=(g_1(v),\ldots,g_d(v))$ with $\exp_{e_0}^{-1}$ gives a $(\nu_i)_{i=1}^d$-dependent coordinate chart 
taking $\bar x := \exp_{e_0} v^r$ to $0 \in \R^d$ and converging in smooth norms as $r \to 0$ to the Riemannian normal coordinates $G^0(v)=(v_1,\ldots,v_d)$ at $e_0$.
 Taylor expansion yields
\[
G^r(v) = DG^r(v^r)(v-v^r) + O(|v-v^r|^2)
\]
where $DG^r(v^r) = I + o(1)$ as $r \to 0$.
Noting inequality \eqref{ell1 model}, given $C' \in (0,1)$,  taking $r \ll s \ll 1$ small enough therefore yields
\begin{align*}
\sum_{i=1}^d |g_i(v)|
& \ge \Big(\sum_{i=1}^d g_i(v)^2\Big)^{1/2}  
\\& \ge \frac{C'+1}2 |v-v^r|    
\\&\ge C' \rho(\bar x,\exp_{e_0} v)\ \  \text{for all} \ \ v \in Q. \nn
\end{align*}
Combining the foregoing, we obtain
\be
\sum_{i=1}^d F_i(v) \ge C\sum_{i=1}^d |g_i(v)| \ge CC'\rho(\bar x,\exp_{e_0} v) \ \  \text{for all} \ \ v \in Q \nn
\ee
which translates to \eqref{agg} on $\bbS^d$ via the exponential map \eqref{spherical exponential}.
\QED

\begin{remark}[Conditions for equality]
Varying $0<C<\frac2\pi$ makes it clear that \eqref{agg} is saturated only when both sides vanish.
\end{remark}

Recall that $\cPb(\bbS^d)$ consists of probability measures $\hat \mu$ whose support covers an orthonormal basis
$V \subset \R^{d+1}$ and is contained in the double $V \cup -V$ of that basis.
For $\al > 1$, 
the following theorem 
provides a $d_\infty$-ball around each such $\hat \mu$ on which it maximizes the energy $E_\al(\mu)$ essentially
uniquely (i.e. uniquely among probability measures on the projective sphere,  apart from rotations). 
It is inspired by 
Corollary~4.3 of \cite{LimMcCann19p}, which gives the analogous result in a different context.

\begin{theorem}[$d_\infty$-local energy maximizers]
\label{localmax}
Given $1<\al <\infty$, $m>0$ and $d \in \bbN$,  there exists $r=r(d,\al, m)>0$ such that for every $\bt \ge \al$
 and $\xi,\hat \xi \in \cP(\bbS^d)$ with $d_\infty(\xi, \hat \xi) < r$:
 if $\hat \xi \in \cPb(\bbS^d)$ and $\hat \xi(\{z,-z\})\ge m$ for each $z\in \spt \hat \xi$,
then $E_{\bt}(\xi) \le E_{\bt}(\hat \xi)$ and the inequality is strict unless $\xi$ is a rotation of $\hat \xi$.
\end{theorem}  

  \noin
{\bf Proof.}  Observe it is sufficient to prove for $\bt = \al$ since $E_\bt(\xi) \le E_\al(\xi)$ for all $\bt\ge \al$ and $\xi \in \cP(\bbS^d)$, while $E_\bt(\hat \xi) = E_\al(\hat \xi)$. Thus we set $\bt=\al$. Fix $\hat \xi \in \cPb(\bbS^d)$ and  assume $\xi \in \cP(\bbS^d)$ satisfies $d_\infty(\hat \xi,\xi) <r$ for some $0<r \ll \pi/4$ to be specified later.
Again the abbreviation $r \ll \pi/4$ means ``for $r$ sufficiently small.''

By rotation, we may assume $\hat \xi = \sum_{i=0}^d ( p_i  \delta_{e_i} + q_i \delta_{-e_i})$. By transferring the mass at $-e_i$ to $e_i$, define $\hat \mu =\sum_{i=0}^d m_i  \delta_{e_i} $ with $m_i = p_i+q_i>0$, 
and set $m:= \min_i m_i$. We similarly transform $\xi$ to $\mu$, retaining $d_\infty(\mu, \hat \mu) < r$ with $E_\al(\hat \xi) = E_\al(\hat \mu)$ and $E_\al(\xi) = E_\al(\mu)$. (We can alternately consider
that we convert $\hat \xi$ and $\xi$ to measures $\hat \omega$ and $\omega$ on the projective space $\bbR\bbP^d$,  
by pushing them forward through the map 
$z \in \bbS^d \to \tilde z=\{z,-z\} \in \bbR\bbP^d$.  This neither increases their separation
nor changes their energy,  when the obvious definitions of $d_\infty$ and $E_\al$ are adopted for 
measures on projective space.  We shall derive conditions under which $\omega$ must be a rotation of $\hat \omega$,
hence supported at $d+1$ well-separated points.  When $d_\infty(\xi,\hat \xi)<\pi/4$ 
this in turn implies $\xi$ is a rotation of $\hat \xi$.)

We need to show that, for $r$ sufficiently small,  $d_\infty(\hat \mu,\mu)<r$ implies $E_\al(\hat \mu) \ge E_\al(\mu)$ and the inequality is strict unless $\mu$ is a rotation of $\hat \mu$, meaning in particular $\mu \in \cPo(\bbS^d)$ as well. 
Note that, since $d_\infty(\mu, \hat \mu) < r \ll \pi/4$, we have $ \mu(D(e_i, r)) = m_i$ for all $i$. Let $\mu_i = \mu \lfloor_{D(e_i, r)}$ be the restriction and $\nu_i=m_i^{-1} \mu_i$ its normalization.
Setting $F(\mu)=2(E_\al(\mu) - E_\al(\hat \mu))$ and $\La^\al * \mu_j(x):= \int \La^\al(x,y) d\mu_j(y)$, observe
\be
F(\mu) = \sum_{j=0}^d \bigg[ \int (\La^\al * \mu_j ) d\mu_j + \sum_{i \neq j}\iint \big( \La^\al(x_i, x_j)-1 \big) d\mu_i(x_i) d\mu_j(x_j) \bigg].   \nn
\ee
By Proposition \ref{P:aggregation}, for any $C \in (0, \frac{2}{\pi})$ we have, recalling $m:=\min_i m_i$,
\begin{align*}
& \sum_{j=0}^d \sum_{i \neq j}\iint \big( \La^\al(x_i, x_j)-1 \big) d\mu_i(x_i) d\mu_j(x_j) \\
 & =  \sum_{j=0}^d \sum_{i \neq j} m_im_j \iint \big( \La^\al(x_i, x_j)-1 \big) d\nu_i(x_i) d\nu_j(x_j) \\
 & \le -Cm^2  \sum_{j=0}^d  \int \dd(\h x_j, x_j) d\nu_j(x_j) 
 \end{align*} 
for some $\{ \h x_j \}_{j=0}^d \subset \bbS^d$, if $r$ is sufficiently small. 

Next, let us address the localized self-interaction terms $\int (\La^\al * \mu_j) d\mu_j$. For $x,x' \in D(e_j, r)$,
\begin{align*}
&\iint \La^\al(x, x') d\mu_j(x)d\mu_j(x') 
 \\&= \Big(\frac {2}{\pi}\Big)^\al m_j^2 \iint \dd(x, x')^\al d\nu_j(x)d\nu_j(x')
\\& \le C' \int \dd(x, \h x_j) ^\al d\nu_j(x)
\end{align*}
with $C' =(\frac{4}{\pi})^\al$, since $\al \ge 1$ yields $ 2^{1-\al} \dd(x, x')^\al \le \dd(x, \h x_j)^\al + \dd(\h x_j, x')^\al$ by convexity. 
Then since $\al >1$ and $\h x_j \to e_j$ as $r \to 0$, for any $\ep >0$ we have
\be
\int \dd(x, \h x_j) ^\al d\nu_j(x) \le \ep \int \dd(x, \h x_j) d\nu_j(x) \nn
\ee
if $r$ is small enough. Combining, we obtain
\be
F(\mu) \le (C'\ep - Cm^2)\sum_{j=0}^d  \int \dd(x, \h x_j) d\nu_j(x). \nn
\ee
By taking $\ep < Cm^2 / C'$ we have $F(\mu) \le 0$  as desired. Moreover $F(\mu) < 0$ unless $\int  \dd(x, \h x_j) d\nu_j(x) =0$ for all $j$, in which case $\nu_j = \delta_{\h x_j}$ for all $j$. Recalling its definition, $F(\mu)=0$ then implies $\dd(\h x_i, \h x_j) =\pi/2$ for all $i \ne j$, that is, $F(\mu)=0$ if and only if $\mu$ is a rotation of $\hat \mu$.  
\QED

\begin{remark}[Sharp threshold exponent for local stability]
When $\al =1$,   the proof of Lemma \ref{iff} shows the conclusion of Theorem \ref{localmax} to fail,  in the sense that even when 
$\hat \xi \in \cPb^=(\bbS^d)$ so that $m=\frac1{d+1}$,
$E_1(\xi)=E_1(\hat \xi)$ does not imply $\xi$ is a rotation of $\hat \xi$,  no matter how small
$d_\infty(\xi,\hat \xi)$ is.
\end{remark}

The remaining sections of this paper establish $\ald<\infty$.  Taking this fact for granted --- at least
temporarily --- yields:

\begin{corollary}[Discontinuous bifurcation unless $\ald=1$]\label{C:bifurcate}
Fix $d \in \bbN$. No curve $(\mu_\al)_{\al>0}$ of optimizers
\begin{equation}\label{optimal curve}
\mu_\al \in \argmax_{\cP(\bbS^d)} E_\al
\end{equation}
can be $d_\infty$-continuous at $\al=\ald$ except possibly if $\ald=1$.
\end{corollary}

\noin 
{\bf Proof.}
Choose any curve $(\mu_\al)_{\al >0}$ of optimizers \eqref{optimal curve}.
Lemma \ref{iff} provides $\ald \ge 1$ such that $\mu_\al \in \cPb^=(\bbS^d)$ if $\al>\ald$,  but not if $\al<\ald$.
Theorem \ref{localmax} provides a $d_\infty$-neigbourhood (in fact, tubular of radius
$r=r(d,\frac{1+\ald}2,\frac1{m+1})>0$) around
$\cPb^=(\bbS^d)$ on which $\cPb^=(\bbS^d)$
uniquely maximizes $E_\al$ for all $\al>  \frac12(1+\ald)$.
Its optimality \eqref{optimal curve} ensures $\mu_\al$ lies outside this neighbourhood whenever $\al<\ald$.
Unless $(\frac{1+\ald}2,\ald)$ is empty, this implies
$$
\liminf_{\epsilon\searrow 0} d_\infty(\mu_{\al+\epsilon},\mu_{\al-\epsilon}) \ge r>0.
$$
Thus $d_\infty$-continuity of $(\mu_\al)_{\al>0}$ at $\al=\ald$ implies $\ald \in \{1,\infty\}$.
But $\ald=\infty$ is ruled out below (see Remark \ref{R:lp}).
\QED

\section{Nonunique energy maximizers at threshold exponent} 
\label{S:nonunique}

In this section we show that the interval of exponents $\al<\infty$ for which the energy \eqref{maximum value} is essentially uniquely
maximized by the conjectured optimizer does not include its endpoint.
We shall employ DeGiorgi's notion of $\Gamma$-convergence~\cite{Braides02}. Since the sign conventions
in this theory are normally set up so that 
$\Gamma$-convergence guarantees accumulation
points of minimizers are minimizers,  we must show $-E_\bt = \Gamma$-$\lim_{\al \to \bt} (-E_\al)$.

\begin{definition}[$\Gamma$-convergence] 
A sequence $F_i:M \longrightarrow \overline\R$ on a metric space $(M,d)$
 is said to {\em $\Gamma$-converge} to $F_\infty:M \longrightarrow \overline\R$,
 denoted $F_\infty  = \Gamma$-$\lim_{i \to\infty}  F_i$, if (a)
\be\label{Gamma-lsc}
F_\infty(\mu) \le \mathop{\lim\inf}\limits_{i\to \infty} F_i(\mu_i) 
\q {\rm whenever} \q d(\mu_i,\mu) \to 0,
\end{equation}
and (b) each $\mu \in M$ is the limit of a sequence $(\mu_i)_{i} \subset M$ along which
\be\label{Gamma-construction}
F_\infty(\mu) \ge \mathop{\lim\sup}\limits_{i\to \infty} {F_i(\mu_i).  }
\end{equation}
\end{definition}

\begin{lemma}[$\Gamma$-convergence of energies]
\label{Gamma1}
 The functionals $(-E_\al)$ {$\Gamma$-converge} to $(-E_{\bt})$ on $(\cP(\bbS^d), d_2)$ as $\al \to \bt \in [1,\infty]$. Here $d_2$ is from \eqref{KRW metric}.
\end{lemma}

\noin 
{\bf Proof.} Let $\{\al_n\}_n$ be a sequence with $\lim_{n \to \infty} \al_n =\bt$. To show the $\Gamma$-convergence, we need to show:
\be\label{Gamma-usc}
E_{\bt}(\mu) \ge \mathop{\limsup}\limits_{n \to \infty} E_{\al_n}(\mu_n) 
\q {\rm whenever} \q d_2(\mu_n,\mu) \to 0,
\end{equation}
and  each $\mu \in \cP(\bbS^d)$ is the limit of a sequence $(\mu_n)_{n} \subset \cP(\bbS^d)$ with
\be\label{Gamma-construction}
E_\bt(\mu) \le \mathop{\liminf}\limits_{n \to \infty} E_{\al_n}(\mu_n). 
\end{equation}
 If $1\le \bt<\infty$ more is true: in this case $E_{\al_n}$ and $-E_{\al_n}$ converge to $E_\bt$ and $-E_\bt$ respectively since whenever $d_2(\mu_n,\mu)\to 0$  then
\begin{align*}
E_{\al_n}(\mu_n) 
&= \int [\La^{\al_n} - \La^{\beta}] d(\mu_n \otimes \mu_n) + \int \La^\bt d(\mu_n \otimes \mu_n)
\\ & \to  \int \La^\bt d(\mu \otimes \mu) =: E_\bt(\mu)
\end{align*}
as $n\to \infty$, noting $\|\La^{\al_n} - \La^\bt \|_\infty \to 0$ and continuity of $\La^\bt$.

By contrast, when $\bt =\infty$, \eqref{Gamma-construction}  follows from $E_{\infty} \le E_{\al_n}$ by taking $\mu_n := \mu$ for all $n$.
Similarly,  the monotone decreasing dependence of the kernel $\La^\al$ on $\al$ implies 
$\limsup_{n \to\infty} E_{\alpha_n}(\mu_n) \le \limsup_n E_{\alpha_m}(\mu_n) = E_{\alpha_m}(\mu)$ for any $m \in \bbN$. Now $m \to \infty$ yields \eqref{Gamma-usc}  by  Lebesgue dominated convergence theorem.   \QED

\begin{remark}[Limits of maximizers maximize limit]
\label{R:limits maximize}
It is well-known \cite{Braides02} (and easy to see) that Lemma \ref{Gamma1} implies that if $\al \to \bt$, any $d_2$-accumulation point of $\mu_\al \in \argmax_{\cP(\bbS^d)} E_\al$ 
must maximize $E_\bt$ on $\cP(\bbS^d)$. 
\end{remark}

\begin{lemma}[Localizing mass of maximizers for approximating problems near that of the limit]\label{EL}
Suppose $\mu_\al$ attains the maximum \eqref{maximum value} 
and $d_2(\mu_\al,\hat \mu) \to 0$ as $\al \to \bt \in (1,\infty)$. 
If $\hat \mu \in \cPb^=(\bbS^d)$ and $r>0$,
then $\spt (\mu_\al) \subset \bigcup_{x \in \spt \hat \mu} (D(x,r) \cup D(-x,r))$ for all $\al$ sufficiently near $\bt$.
 \end{lemma}

\noin{\bf Proof.} 
Let $\{v_i\}_{i=0}^d$ be the orthonormal basis of $\R^{d+1}$ such that $\spt(\hat \mu) \subset \bigcup_{i=0}^d\{v_i, -v_i\}$.
We first claim the potential function $(\La^\bt * \hat \mu)(x) = \int \La^\bt(x,y) d\mu(y)$ attains its maximum precisely on the set $\bigcup_{i=0}^d\{v_i, -v_i\}$.
Indeed, $\La^\bt *\hat \mu$ is an equally weighted sum of $d+1$ potentials $\La^\bt * \delta_{v_i}$,  each of which is smooth and has a (non-vanishing) non-negative definite Hessian 
away from ($\pm v_i$ and) the equatorial great sphere
$v_i^\perp \cap \bbS^d$ --- which coincides with the cut locus of $v_i$ on $\bbP\bbR^d$.
Thus the only local maxima of $\La^\bt *\hat \mu$ on $\bbS^d$ must lie on $\cup_{i=0}^d v_i^\perp$.
However $\La^\bt * \delta_{v_0}$ attains its maximum value of $1$ throughout the great sphere
 $\bbS^d \cap v_0^\perp$, and it is easy to see by induction on dimension that the local maxima of 
$\sum_{i=1}^d \La^\bt * \delta_{v_i}$ on $\bbS^d \cap v_0^\perp$ occur precisely at $v_1,\ldots,v_d$, to establish the claim.

Since the kernels $\La^\al$ are equiLipschitz for $\al$ near $\bt$, 
as $\al \to \bt$, the potential functions $\La^\al * \mu_\al$ converge uniformly to $\La^\bt *\hat \mu $ on $\bbS^d$. Given $r>0$, this implies that the maximum of $\La^\al * \mu_\al$ occurs only in $\bigcup_{i=0}^d(D(v_i,r) \cup D(-v_i,r))$ for all $\al$ sufficiently near $\bt$. Then  the Euler-Lagrange equation due to the fact $\mu_\al  \in \argmax E_\al$ asserts
\be
\spt (\mu_\al) \subset \argmax \La^\al * \mu_\al, \nn
\ee 
yielding the lemma. 
\QED

\begin{remark}[The mildest repulsion limit]\label{R:mildest exception}
Although we expect Lemma \ref{EL} to remain true when $\bt=\infty$,  our proof does {not} extend to that case since the kernels $\Lambda^\al$ 
fail to be equiLipschitz and the continuous potentials $\Lambda^\al * \mu_\al$ cannot converge uniformly to the discontinuous limit 
$\Lambda^\infty *\hat \mu$ as $\al \to \infty$.
Notice however, that if Lemma \ref{EL} holds for $\bt=\infty$, then the next proof shows Theorem \ref{main1} also extends to $\ald=\infty$.
Such an extension will be relevant in Section \ref{S:discrete}.
\end{remark}

\begin{theorem}[Maximizers are not unique at threshold exponent]
\label{main1}
If $\ald < \infty$, then $\cPb^=(\bbS^d) \subset \argmax_{\cP(\bbS^d)} E_\ald$, and the inclusion is strict.
\end{theorem}

\noin{\bf Proof.} 
For $\ald < \infty$,  the desired inclusion \eqref{nonunique max} follows from the definition of $\ald \ge 1$ in Lemma \ref{iff} using Remark \ref{R:limits maximize}.
Our goal is therefore to show this inclusion is strict.
To derive a contradiction,  suppose $\cPb^=(\bbS^d) = \argmax_{\cP(\bbS^d)} E_\ald$. 
The compactness of $\cP(\bbS^d)$ and continuity of $E_\al$ (for the metric $d_2$),
provide $\mu_\al \in \argmax_{\cP(\bbS^d)} E_\al$ for each $\al\in(0,\infty)$. 
The same compactness provides an increasing
 sequence $\al_n \nearrow \ald$ and corresponding maximizers $\mu_{\al_n} \in \argmax E_{\al_n}$ which $d_2$-converge to some $\hat \mu \in \argmax E_\ald$. 
By assumption $\hat \mu \in\cPb^=(\bbS^d)$. Then the Localization Lemma \ref{EL} and Local Stability Theorem \ref{localmax} combine to yield that for large enough $n$,
 there exists $\nu_{\al_n} \in \cPb(\bbS^d)$ obtained by collapsing each bit of the mass of $\mu_{\al_n}$ onto the nearest point in the support of $\hat \mu$ such that {$\spt (\nu_{\al_n}) = \spt(\hat\mu)$,} $d_\infty(\mu_{\al_n}, \nu_{\al_n})$ is small and
 $E_{\al_n} (\mu_{\al_n}) \le E_{\al_n}(\nu_{\al_n})$, with the inequality being strict unless $\mu_{\al_n} \equiv \nu_{\al_n}$. 
Moreover, {since $\spt (\nu_{\al_n}) = \spt(\hat\mu)$,} {it is well-known that}
 $E_{\al_n}(\nu_{\al_n}) \le E_{\al_n}(\hat \mu)$ and the inequality is strict unless $\nu_{\al_n} \equiv\hat \mu$, i.e. $\nu_{\al_n} \in \cPb^=(\bbS^d)$, by e.g. Perron-Frobenius \eqref{Perron-Frobenius}. 
 Since $\mu_{\al_n} \in \argmax E_{\al_n}$, the inequalities must be equalities, hence $\cPb^=(\bbS^d) \subset \argmax_{\cP(\bbS^d)} E_{\al_n}$ for all large $n$. This implies $\ald \le \al_n$ by Lemma \ref{iff}, contradicting the strict monotonicity $\al_n \nearrow \ald$ to establish the theorem. 
 \QED

\begin{remark}[Finiteness of the threshold exponent]\label{R:lp} 
In the first version of this manuscript posted on the arXiv,  we combined similar arguments with our prior results \cite{LM20_1} \cite{LM20_3} 
and a uniform cardinality bound on the support of $d_2$-local maximizers for all $\al \ge 4$ to deduce $\ald <\infty$; see the next section for a similar but simpler argument in the discrete context $N<\infty$.
Thus  the finiteness required by Theorem \ref{main1} turns out not to be restrictive.  
In response, Bilyk et al 
obtained an explicit bound $\ald \le 2$ using a clever majorization argument,
which they were kind enough to include in the attached appendix.
Combining their work with ours yields a stronger conclusion,
which might be seen as additional evidence for the continuous Fejes-T\'oth conjecture $\ald=1$, namely:
\end{remark}

\begin{corollary}[Improving Bilyk et al's threshold exponent bound]\label{C:less than 2}
The threshold exponent from Lemma \ref{iff} satisfies $\ald \in [1,2)$.
\end{corollary}

\noin{\bf Proof.} Follows from Theorems \ref{main1} and \ref{T:majorization} and Lemma \ref{iff}.
\QED
\\

Finally,  we remark on the possible disconnectedness of the set of optimizers at the threshold exponent.
Although the set $\cPb^=(\bbS^d)$ of measures on the sphere has uncountably many $d_\infty$-connected components,  they are all essentially equivalent to 
a single measure $\hat \mu$.  
The rotates of $\hat \mu$ form a $d_\infty$-connected set of measures 
 on projective space $\bbR \bbP^d = \bbS^d/\{\pm\}$. 
 However,   the results of this section and the previous one 
combine to show $d_\infty$-connectedness of the set of energy optimizers on projective space must fail 
at the threshold exponent unless the continuous Fejes T\'oth conjecture holds true:

\begin{corollary}[On connectedness of the threshold set of optimizers]\label{C:disconnect}
Fix $d \in \bbN$. At the threshold exponent $\ald$, 
 if the set of optimizers 
\be\label{threshold optimizers}
(\argmax_{\cP(\bbS^d)} E_\ald) / \{\pm\}
\ee
 forms a $d_\infty$-connected subset of $\cP(\bbR \bbP^d)$, then $\ald=1$.
\end{corollary}

\noin{\bf Proof.} Assume $\ald>1$.  Theorem \ref{localmax} then provides a $d_\infty$-open neighbourhood of $\cPb^=(\bbS^d)$ on
which $\cPb^=(\bbS^d)$ uniquely maximizes $E_\ald$.  Theorem \ref{main1} and the finiteness of $\ald$ from Remark \ref{R:lp} assert another maximizer
lies outside this neighbourhood.  The $d_\infty$-disconnectedness of \eqref{threshold optimizers} follows,
and establishes the corollary's contrapositive.
\QED
\\

Of course, if the conjecture $\ald=1$ holds true, we do not know whether or not the set of optimizers \eqref{threshold optimizers} is 
$d_\infty$-connected.
Indeed, the last corollary 
shows the continuous Fejes T\'oth conjecture follows if one can establish 
$d_\infty$-connectedness of \eqref{threshold optimizers}.

\section{Maximizing energy among finitely many lines (charges)}\label{S:discrete}

 Fejes T\'oth's original paper \cite{FT59} considered the optimization problem $\max E_1$ not on $\cP(\bbS^d)$, but instead
 over the domain of uniform discrete probabilities $\cP_N^=(\bbS^d)$ for each $N \in \bbN$. It conjectured
\begin{align*}
&\cP_{N,\Delta}^=(\bbS^d):=\{ \mu \ | \ \text{there is an orthonormal basis $\{v_i\}_{i=1}^{d+1}$ of $\R^{d+1}$ and }  \\
&\{x_i\}_{i=1}^N \text{ such that $ \mu = \frac{1}{N}\sum_{i=1}^N \delta_{x_i}$ and $x_i \in \{v_j, -v_j\}$ if $i \equiv j$ mod ${d+1}$} \} 
\end{align*}
maximizes $E_1$ over $\cP_N^=(\bbS^d)$. Notice $\mu \in \cP_{N,\Delta}^=(\bbS^d)$ iff $\mu \equiv \psi$ where $\psi$ distributes $N$ particles of equal mass $1/N$ on an orthonormal basis as uniformly as possible.

This problem admits an analog of Lemma \ref{iff},  whose similar proof is omitted:

\begin{lemma}[Threshold exponent for discrete equidistribution over maximal simplices]\label{L:N-iff} 
Given $d,N \in \bbN$, there exists a unique $\ald(N) \in [0,\infty]$ 
such that $0< \al < \ald(N)$ implies $\cP_{N,\Delta}^=(\bbS^d)$ disjoint from $\argmax_{\cP_N^=(\bbS^d)} E_\al$, whereas $\ald<\al$ implies
\be\label{N-unique max}
\cP_{N,\Delta}^=(\bbS^d) = \argmax_{\cP_N^=(\bbS^d)} E_\al.
\ee
If $N >d+1$ then \eqref{N-unique max} fails for $\al=1$, hence  $\ald(N) \ge 1$.
\end{lemma}

In terms of this $N$-dependent threshold exponent,  Fejes T\'oth's conjecture reduces to the assertion that $N>d+1$ implies $\ald(N)=1$.

The majorization method of the appendix below 
does not seem to apply in this setup unless $N$ is a multiple of $d+1$. 
Thus for $N$ indivisible by $d+1$ it is not clear to us whether  $\ald(N) < 2$. However, we shall adapt the strategy above to show that $\ald(N) < \infty$, that is for each $N$, 
\eqref{N-unique max} holds for all sufficiently large finite $\al$.

To this end, we first recall our results \cite[Theorem 1.4]{LM20_1} on the mildest repulsion limit $\al = \infty$, which include the assertion that \eqref{N-unique max} holds when $\al=\infty$. 
{(We subsequently learned from the authors of our appendix that this conclusion can alternately be obtained as a consequence of Turan's theorem, which asserts that, among graphs with $N$ vertices and no $(d+2)$-cliques, the number of 
edges is uniquely maximized by
the complete $(d+1)$-partite graph
$K_{k_1,\ldots,k_{d+1}}$ with $|k_i-k_j|\le 1$ for all $i,j \in \{1,\ldots, d+1\}
$
\cite{Aigner95} \cite{Turan41}.)}

The proof of the following lemma is virtually verbatim to that of Lemma~\ref{Gamma1}, 
if one replaces $E_\bt$ with $E_\infty$ and $\cP(\bbS^d)$ with $\cP_N^=(\bbS^d)$, so we omit it here.

\begin{lemma}[$\Gamma$-convergence of discrete energies]
\label{Gamma2}
Fix $N \in \bbN$. The functionals $(-E_\al)$ {$\Gamma$-converge} to $(-E_\infty)$ on $(\cP_N^=(\bbS^d),d_2)$ as $\al \to \infty$. Here $d_2$ is from \eqref{KRW metric}.
\end{lemma}

As in Remark \ref{R:limits maximize}, this implies any $d_2$-accumulation point $\mu$ of $\mu_\al \in \argmax_{\cP_N^=(\bbS^d)} E_\al$ as $\al \to \infty$  belongs to $\argmax_{\cP_N^=(\bbS^d)} E_\infty = \cP_{N,\Delta}^=(\bbS^d)$. 
 But since $\mu_\al \in \cP_N^=(\bbS^d)$, the following lemma (which may be folklore to experts)
 shows $d_2(\mu_\al, \mu) \to 0$ yields $d_\infty(\mu_\al, \mu) \to 0$ as $\al \to \infty$.

\begin{remark}
The sphere $\bbS^d$ plays no special role,  and could be replaced by an arbitrary complete separable metric space in the next lemma.
\end{remark}

\begin{lemma}[$d_\infty$ metrizes the narrow topology on $\cP_N^=(\bbS^d)$]
\label{L:discrete narrow}
Given $1 \le p <\infty$, any $\mu,\nu \in \cP_N^=(\bbS^d)$ satisfy
\[
N^{-1/p} d_\infty(\mu,\nu) \le d_p(\mu,\nu) \le d_\infty(\mu,\nu).
\]
\end{lemma}

\noin{\bf Proof.}
One direction is standard:
Jensen's inequality implies $d_p \le d_q$ for all $1\le p< q <\infty$,  hence $d_p \le d_\infty = \lim_{q \to \infty} d_q$
even on the larger space $\cP(\bbS^d)$.  Conversely,  
$\mu,\nu \in \cP_N^=(\bbS^d)$ satisfy $\mu = \frac 1N \sum_{i=1}^N \delta_{x^i}$
and $\nu = \frac 1N \sum_{j=1}^N \delta_{y^j}$, so
 the infimum \eqref{KRW metric} over $\gamma \in \Gamma(\mu,\nu)$ 
defining $d_p(\mu,\nu)$ reduces to a finite-dimensional linear program
for a doubly stochastic $N \times N$ matrix with entries $G^{ij}=\gamma[\{(x^i,y^j)\}]$.
This infimum is attained by a permutation matrix $G_p$,  as a consequence of the
Birkhoff-von Neumann characterization of extremal doubly stochastic matrices 
\cite{Birkhoff46}.
Reorder the points $(y^j)_{j=1}^N$ so that $G_p$ is the identity matrix.
Taking the corresponding measure $\gamma_p\in\Gamma(\mu,\nu)$ as a competitor bounding the distance $d_q(\mu,\nu)$,  
the limit $q \to \infty$  (or \eqref{sand metric} directly) yields
\begin{align*}
N d_p(\mu,\nu)^p 
&= \sum_{i=1}^N \rho(x^i,y^i)^p   
\\&\ge  \max_{1 \le i \le N} \rho(x^i,y^i)^p
\\&\ge d_\infty(\mu,\nu)^p
\end{align*}
as desired.
\QED
\\

Using this $d_\infty$-convergence in place of Lemma \ref{EL} and noting Remark~\ref{R:mildest exception}, the same proof adapts
Theorem \ref{localmax} to this context, yielding the following conclusion:

\begin{theorem}[Finiteness of discrete threshold exponent]
\label{main2}
Given $d, N \in \bbN$, the threshold exponent from Lemma \ref{L:N-iff} is finite $\ald(N) \in [0,\infty)$ and the 
containment 
$\cP_{N,\Delta}^=(\bbS^d) \subset \argmax_{\cP_N^=(\bbS^d)} E_{\ald(N)}$ is strict.
\end{theorem}

{The obvious analogs of the discontinuity / disconnectedness results of
Corollaries \ref{C:bifurcate} and \ref{C:disconnect} extend to this discrete (i.e.~$N<\infty$) setting.}

\appendix
\section{
by Dmitriy Bilyk, Alexey Glazyrin, Ryan Matzke, Josiah Park, and Oleksandr
 Vlasiuk}

\begin{theorem}[Estimating the threshold exponent: $\ald \le 2$] \label{T:majorization}
For every $\al \ge 2$, the set of maximizers of \eqref{maximum value} is precisely $\cPb^=(\bbS^d)$.
\end{theorem}

\noin{\bf Proof.}  {For $t \in [-1,1]$ set $f_\al (t) = (\frac{2}{\pi} \arccos |t| )^\al$ and $g(t)=1-t^2$. We claim
\be\label{majorize} 
g \ge f_\al \text{ on } [-1,1] \text{ and } \{g=f_\al\} =  \{-1,0,1\} \iff \al \ge 2.
\ee
Indeed, setting $h_\al(t) := g(t)^{1/\al} - f_\al(t)^{1/\al}$,  the computation $h_2''(t) = (2t\pi^{-1} -1)(1-t^2)^{-3/2}$ shows $h_2(t)$ to be
strictly concave on the interval $t \in [0,1]$ and to vanish at both endpoints.  For  $\al = 2$ this establishes \eqref{majorize}.
For $\al>2$, the same conclusion then follows from $\frac{dh_\al}{d\al}(t)\ge 0$.
For $\al<2$ we have {$\lim_{t \nearrow 1} f_\al'(t)=- \infty<g'(1)$ and $f_\al(1)=0=g(1)$, thus domination of $f_\al$ by $g$ fails,} 
confirming the reverse implication in \eqref{majorize}.

Now \eqref{energy} may be rewritten in the form $E_\al(\mu) = {F_{f_\al}}(\mu)$ where
\[
{F_f}(\mu) := \iint  f(x \cdot y) d\mu(x)d\mu(y).
\]
Let $\mu \in \cP(\bbS^d)$, $\sigma$ be the uniform probability on $\bbS^d$, and $\hat \mu \in \cPb^=(\bbS^d)$. 
For $\al \ge 2$ we claim
\be\label{Bilyk chain}
{F_{f_\al}}(\mu) \le {F_{g}}(\mu) \le {F_g}(\sigma) = {F_g}(\hat \mu) = {F_{f_\al}} (\hat \mu) 
\ee
where the middle two (in)equalities for $g$ are known and reproved below, while the first and last follow from \eqref{majorize}
--- which also makes the first inequality strict unless $\mu$ lies in the narrow closure of $\cPb(\bbS^d)$. 
On the other hand, for $\mu \in \overline{\cPb(\bbS^d)}$, \eqref{Perron-Frobenius} implies that
the second inequality in \eqref{Bilyk chain}
becomes strict unless $\mu \in \cPb^=(\bbS^d)$. 

It remains to establish
the middle two (in)equalities in \eqref{Bilyk chain} for $g$, which can be done in various ways, c.f.~\cite{BenedettoFickus03}.
For example,  defining the symmetric $(d+1) \times (d+1)$ matrix $I(\mu)$ by
\[
I^{ij}(\mu) = \int_{\bbS^d} x^i x^j d\mu(x),
\]
the Cauchy-Schwartz inequality for the Hilbert-Schmidt norm yields
\begin{align*}
\iint (x \cdot y)^2 d\mu(x)d\mu(y) 
&= \Tr(I(\mu)^2) 
\\&\ge \frac{1}{d+1} (\Tr I(\mu))^2
\\&= \frac1{d+1}
\end{align*}
since $\Tr I(\mu) = |x|^2 = 1$ on $\bbS^d$.  If $I(\mu)$ is a multiple of the identity matrix --- as for $\mu\in \{\hat \mu,\sigma\}$ --- then equality holds.  This establishes \eqref{Bilyk chain} and completes the proof.
\QED

\section*{Declarations}

The authors are not aware of any conflicts of interest;  
they have no financial or proprietary interests in any material discussed in this article.
TL is a faculty member in Purdue University's Krannert School of Management.  RJM is a Canada Research Chair at the University of Toronto.  TL's work was also supported in part by ShanghaiTech University,  the University of Toronto and its Fields Institute for the Mathematical Sciences.  RM acknowledges partial support of this research by the Canada Research Chairs Program and Natural Sciences and Engineering Research Council of Canada Grants RGPIN 2015-04383 and 2020-04162. 
Data sharing is not applicable to this article as no datasets were generated or analysed during the current study.

\end{document}